\documentclass[a4paper,10pt]{article}
\usepackage[utf8x]{inputenc}
\usepackage{tracefnt,amsmath,tabu,array}
\usepackage{amssymb,graphicx,setspace,amsfonts,amsbsy}
\usepackage{pifont,latexsym,ifthen,amsthm,rotating,calc,textcase,booktabs,cancel,slashed}

\newtheorem{theorem}{Theorem}[section]
\newtheorem{lemma}[theorem]{Lemma}

\newtheorem{proposition}[theorem]{Proposition}

\newtheorem{remark}[theorem]{Remark}
\newcommand{\filledbox}{\leavevmode
  \hbox to.77778em{%
  \hfil\vbox to.675em{\hrule width.6em height.6em}\hfil}}

\newcommand{\Rm}{{\mathbb R}}

\begin{document}
\tabulinesep=1.0mm
\title{Channel-localized Strichartz estimates of radial wave equations}

\author{Liang Li, Shenghao Luo and Ruipeng Shen\\
Centre for Applied Mathematics\\
Tianjin University\\
Tianjin, China
}

\maketitle

\begin{abstract}
 In this work we give a few new Strichartz estimates of radial solutions to the wave equation. These Strichartz estimates still use $L^p L^q$ type norms in each channel-like region $\{(x,t): |t|+2^k < |x| < |t|+2^{k+1}\}$, with weaker restrictions on $p, q$ than the classic ones, but combine these localized norms together in the way of an $l^2$ space. We also give an application of these Strichartz estimates on the well-posedness theory of non-linear wave equations.
\end{abstract}

\section{Introduction}

\paragraph{Background} Strichartz-type estimates play an essential role in the well-posedness theory of non-linear wave equations. The following version of Strichartz estimates given by Ginibre-Velo \cite{strichartz} are frequently used in the discussion of non-linear wave equations. 

\begin{proposition}[Strichartz estimates, here we use the Sobolev version]\label{Strichartz estimates} 
Assume $d\geq 2$. Let $2\leq p_1,p_2 \leq \infty$, $2\leq q_1,q_2 < \infty$ and $\rho_1,\rho_2,s\in \Rm$ be constants with
 \begin{align*}
  &\frac{2}{p_i} + \frac{d-1}{q_i} \leq \frac{d-1}{2},& &(p_i,q_i)\neq \left(2,\frac{2(d-1)}{d-3}\right),& &i=1,2;& \\
  &\frac{1}{p_1} + \frac{d}{q_1} = \frac{d}{2} + \rho_1 - s;& &\frac{1}{p_2} + \frac{d}{q_2} = \frac{d-2}{2} + \rho_2 +s.&
 \end{align*}
 Assume that $u$ is the solution to the linear wave equation
\[
 \left\{\begin{array}{ll} \partial_t u - \Delta u = F(x,t), & (x,t) \in \Rm^d \times [0,T];\\
 (u,u_t)|_{t=0} = (u_0,u_1) \in \dot{H}^s(\Rm^d)\times \dot{H}^{s-1} (\Rm^d). & 
 \end{array}\right.
\]
Then we have
\begin{align}
 \left\|\left(u(\cdot,T), \partial_t u(\cdot,T)\right)\right\|_{\dot{H}^s \times \dot{H}^{s-1}} & +\|D_x^{\rho_1} u\|_{L^{p_1} L^{q_1}([0,T]\times \Rm^d)} \nonumber\\
 & \leq C\left(\left\|(u_0,u_1)\right\|_{\dot{H}^s \times \dot{H}^{s-1}} + \left\|D_x^{-\rho_2} F(x,t) \right\|_{L^{\bar{p}_2} L^{\bar{q}_2} ([0,T]\times \Rm^d)}\right). \label{inequality general Strichartz}
\end{align}
Here the coefficients $\bar{p}_2$ and $\bar{q}_2$ satisfy $1/p_2 + 1/\bar{p}_2 = 1$, $1/q_2 + 1/\bar{q}_2 = 1$. The constant $C$ does not depend on $T$ or $u$. 
\end{proposition}

\noindent Let us consider the conditions that $(p_i,q_i)$ must satisfy. The conditions consist of two parts

\begin{figure}[h]
 \centering
 \includegraphics[scale=1.25]{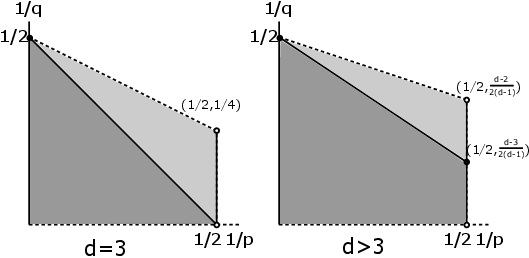}
 \caption{Admissible regions} \label{figure admissible}
\end{figure}

\paragraph{Rescaling identity} The identities 
\begin{align*}
 &\frac{1}{p_1} + \frac{d}{q_1} = \frac{d}{2} + \rho_1 - s;& &\frac{1}{p_2} + \frac{d}{q_2} = \frac{d-2}{2} + \rho_2 +s.&
\end{align*}
are actually determined by the natural rescaling of the wave equation. More precisely, given any linear free wave $u(x,t)$ with initial data $(u_0,u_1)\in \dot{H}^s (\Rm^d) \times \dot{H}^{s-1}(\Rm^d)$ and any positive constant $\lambda > 0$, the function $u_\lambda(x,t) = \lambda^{-d/2+s} u(x/\lambda, t/\lambda)$ is also a linear free wave with initial data $(\lambda^{-d/2+s} u_0(x/\lambda), \lambda^{-d/2-1+s}u_1(x/\lambda))$. These two pair of initial data share the same $\dot{H}^s(\Rm^d)\times \dot{H}^{s-1}(\Rm^d)$ norm. A straight-forward calculation shows that 
\[
 \|D_x^{\rho_1} u_\lambda\|_{L^{p_1} L^{q_1}(\Rm \times \Rm^d)} = \lambda^{-d/2+s-\rho_1+1/p_1+d/q_1} \|D_x^{\rho_1} u\|_{L^{p_1} L^{q_1}(\Rm \times \Rm^d)}.
\]
If the Strichartz estimates \eqref{inequality general Strichartz} hold, then the invariance of $\dot{H}^s \times \dot{H}^{s-1}$ norm implies the uniform boundedness of $\|D_x^{\rho_1} u_\lambda\|_{L^{p_1} L^{q_1}(\Rm \times \Rm^d)}$ for all positive constants $\lambda > 0$. This immediately gives the identity
\[  
 -d/2+s-\rho_1+1/p_1+d/q_1 = 0.
\]
The proof of identity $(d-2)/2 + \rho_2 +s - 1/p_2 - d/q_2 = 0$ is similar if we incorporate the rescaled inhomogeneous term $F_\lambda = \lambda^{-d/2-2+s} F(x/\lambda, t/\lambda)$.   
\paragraph{Admissible region} We fix a dimension $d\geq 3$. The pairs $(p_i,q_i)$ also need to satisfy a few inequalities. This is equivalent to saying that $(1/p_i,1/q_i)$ is contained in the following region 
\begin{align*}
 \left\{\left(\frac{1}{p}, \frac{1}{q}\right) \in \left[0,\frac{1}{2}\right]\times \left(0,\frac{1}{2}\right]: \frac{2}{p}+\frac{d-1}{q} \leq \frac{d-1}{2}; \left(\frac{1}{p}, \frac{1}{q}\right) \neq \left(2,\frac{2(d-1)}{d-3}\right) \right\}.
\end{align*}
We usually call it the admissible region. The restriction $\left(\frac{1}{p}, \frac{1}{q}\right) \neq \left(2,\frac{2(d-1)}{d-3}\right)$ can be further removed in dimension $d\geq 4$ (but not for dimension $3$), according to Keel-Tao \cite{endpointStrichartz}. The admissible region can be further extended if we consider special free waves. For example, if the free waves are radially symmetric in the spatial variable, then the Strichartz estimates 
\[
 \|u\|_{L^{p} L^{q}(\Rm\times \Rm^d)}  \leq C \left\|(u_0,u_1)\right\|_{\dot{H}^s \times \dot{H}^{s-1}} 
\]
holds if $(p,q)$ satisfies the rescaling identity $1/p+d/q = d/2-s$ and is contained in the admissible region 
\[
 \left\{\left(\frac{1}{p}, \frac{1}{q}\right) \in \left[0,\frac{1}{2}\right]\times \left(0,\frac{1}{2}\right]: \frac{1}{p}+\frac{d-1}{q} \leq \frac{d-1}{2}\right\}.
\]
Please see Sterbenz \cite{angularStrichartz}, for example. The admissible region can be visually illustrated as in figure \ref{figure admissible}, where the darker region is the admissible region for general free waves, the lighter region is the additional part for radial free waves. 

\paragraph{Motivation} Let us consider the local theory, especially the continuous dependence of solutions on the initial data, of the energy-critical wave equation 
\[
 \partial_t^2 u - \Delta u = \pm |u|^{4/(d-2)} u.
\]
Basic Fourier analysis gives the Strichartz estimates
\[
 \sup_{t} \|(u, u_t)\|_{\dot{H}^s \times \dot{H}^{s-1}} \leq \|(u_0,u_1)\|_{\dot{H}^s \times \dot{H}^{s-1}} + \|F\|_{L^1 L^2}.
\] 
The Strichartz norm $\|u\|_{L^p L^q}$ with $(p,q) = (\frac{d+2}{d-2}, \frac{2(d+2)}{d-2})$ seems to be a natural choice in the well-posedness theory since we always have (for convenience we use the notation $F(u) = \pm |u|^{4/(d-2)} u$)
\begin{align*}
 &\|F(u)\|_{L^1 L^2} \leq \|u\|_{L^p L^q}^\frac{d+2}{d-2};&
 &\|F(u)-F(v)\|_{L^1 L^2} \leq C (\|u\|_{L^p L^q}^\frac{4}{d-2} +\|v\|_{L^p L^q}^{\frac{4}{d-2}}) \|u-v\|_{L^p L^q}.&
\end{align*}
This Strichartz norm works well in the lower dimensional case $3 \leq d \leq 6$. Proposition \ref{Strichartz estimates} immediately gives 
\[
 \|u\|_{L^p L^q} \lesssim \|(u_0,u_1)\|_{\dot{H}^1 \times L^2} + \|F\|_{L^1 L^2}. 
\]
A standard fixed-point argument then gives the existence and uniqueness of solutions. We may also apply a boot stripe argument to obtain the continuous dependence of solutions on the initial data. In higher-dimensional case, however, the pair $(p,q)$ given above is no longer contained in the admissible region, even in the radial setting, since we have $p<2$ in this case. As a result, it is common to utilize the Strichartz estimates with fractional derivative in higher dimensions. Please see, Duyckaerts-Kenig-Merle \cite{oddtool}, for example. The fractional derivative often leads to technical difficulties when the properties of solutions are discussed. For example, in order to show the continuous dependence of solutions on the initial data, it is necessary to give an upper bound of $\|D^\alpha (F(u)-F(v))\|$ in terms of $\|D^\alpha (u-v)\|$ and $\|u-v\|$. Here $F(u) = |u|^{4/(d-2)} u$. If $\alpha\in (0,1)$ and $d\leq 6$, it is known that 
\begin{align*}
 \|D^\alpha (F(u)-F(v))\|_{L^{p}} \lesssim & (\|F'(u)\|_{L^{p_1}} + \|F'(v)\|_{L^{p_1}}) \|D^\alpha(u-v)\|_{L^{p_2}}  \\
 & + (\|F''(u)\|_{L^{r_1}} + \|F''(v)\|_{L^{r_1}}) (\|D^\alpha u\|_{L^{r_2}} + \|D^\alpha v\|_{L^{r_2}}) \|u-v\|_{L^{r_3}}. 
\end{align*}
Here $1/p = 1/p_1+1/p_2 = 1/r_1 + 1/r_2 + 1/r_3$. Please see Kenig-Merle \cite{kenig} and citations therein. If $d\geq 7$, however, this does not hold because of the singularity of $F''$ near zero. In addition, the fractional derivative is also inconvenient when we apply a spatial cut-off operator. This kind of operator is frequently used in the discussion of wave equations, due to the finite speed of wave propagation. In summary, it is attractive to use a Strichartz norm which is a slightly modified version of the $L^p L^q$ norm. Inspired by the successful channel of energy method, we first consider the $L^p L^q$ type norm localized in a channel-like region $\{(x,t): |t| + 2^k < |x| < |t|+ 2^{k+1}\}$ and then combine them together in the way of an $l^2$ space. Please note that our new Strichartz estimates concern the values of solutions outside the main light cone $|x|=|t|$ thus are especially useful in the discussion of asymptotic behaviours of solutions near the infinity. 

\paragraph{Notations} Before we give the statement of our main result, we first introduce a few notations for convenience. We define the regions
\[
 \Omega_j = \{(x,t)\in \Rm^d \times \Rm: |t|+2^j \leq |x| < |t|+2^{j+1}\}. 
\]
These regions are disjoint and their union is exactly the exterior region $\{(x,t): |x|>|t|\}$ outside the main light cone $|x|=|t|$. We also use the notation $\chi_j$ for the characteristic function of $\Omega_j$. We define the space-time norm ($I\subset \Rm$ is an arbitrary time interval)
\[
 \|u\|_{L^p L^q (I \times \Rm^d)} = \left(\int_{I} \left(\int_{\Rm^d} |u(x,t)|^p dx \right)^{q/p}dt\right)^{1/q}.
\]

\paragraph{Main result} Now we give our main result of this work and a few remarks. Our Strichartz estimates use modified $L^p L^q$ type norms with very weak assumptions on $p,q$. 

\begin{theorem} 
Assume that $d\geq 3$ is odd. Let $u$ be a radial solution to the linear wave equation 
\[
 \left\{\begin{array}{ll} \partial_t^2 u - \Delta u = F, & (x,t) \in \Rm^d \times \Rm; \\ (u, u_t)|_{t=0} = (u_0,u_1). & \end{array} \right.
\]
In addition, we assume the constants $\beta \in (1/2, 3/2)$ and $p,q, \tilde{p}, \tilde{q} \in [1,+\infty)$ satisfy 
\begin{align*}
 &\frac{1}{p} + \frac{d}{q} = \frac{d}{2} - \beta; & &\frac{1}{\tilde{p}} + \frac{d}{\tilde{q}} = \frac{d}{2} - \beta + 2.&
\end{align*}
Then the following inequality holds
\[
 \left(\sum_{j=-\infty}^\infty \|\chi_j u\|_{L^p L^q (\Rm \times \Rm^d)}^2 \right)^{1/2} \leq C \left[ \|(u_0,u_1)\|_{\dot{H}^\beta \times \dot{H}^{\beta-1}(\Rm^d)} + \left(\sum_{k=-\infty}^\infty \|\chi_k F\|_{L^{\tilde{p}} L^{\tilde{q}} (\Rm \times \Rm^d)}^2 \right)^{1/2} \right].
\]
The constant $C$ above does not depend on the initial data $(u_0,u_1)$ or the inhomogeneous term $F$. 
\end{theorem}

\begin{remark}
By finite speed of propagation we also have for any $N \in \mathbb{Z}$: 
 \[
 \left(\sum_{j=N}^\infty \|\chi_j u\|_{L^p L^q (\Rm \times \Rm^d)}^2 \right)^{1/2} \leq C \left[ \|(u_0,u_1)\|_{\dot{H}^\beta \times \dot{H}^{\beta-1}(\Rm^d)} + \left(\sum_{k=N}^\infty \|\chi_k F\|_{L^{\tilde{p}} L^{\tilde{q}} (\Rm \times \Rm^d)}^2 \right)^{1/2} \right].
\]
In addition, we may substitute $\Rm$ by any time interval $I$ containing zero. 
\end{remark}

\begin{remark}
 The Strichartz estimates given in the main theorem concern the values of radial solution $u$ in the exterior region $\{(x,t): |x|>|t|\}$. This region touches the axis $x=0$ at a single point. In fact, if $\Omega\subset \Rm^d \times \Rm$ is a region containing an open neighbourhood of $(0,t_0)$ and the Strichartz estimates ($\chi_\Omega$ is the characteristic function of $\Omega$, $p+d/q = d/2 -\beta$)
 \[
  \|\chi_\Omega u\|_{L^p L^q (\Rm \times \Rm^d)} \lesssim \|(u_0,u_1)\|_{\dot{H}^\beta \times \dot{H}^{\beta-1}(\Rm^d)}
 \]
 holds for all radial free waves $u = \mathbf{S}_L(u_0,u_1)$, then we may apply a time transformation and a rescaling to conclude that the global Strichartz estiamte
 \[
  \|u\|_{L^p L^q (\Rm \times \Rm^d)} \lesssim \|(u_0,u_1)\|_{\dot{H}^\beta \times \dot{H}^{\beta-1}(\Rm^d)}
 \]
 holds for all radial free waves. 
\end{remark}

\begin{remark}
 Let us consider the pair $(p,q) = \left(\tfrac{d+2}{d-2}, \tfrac{2(d+2)}{d-2}\right)$ mentioned above. If $3\leq d \leq 6$, then the Strichartz estimates 
 \[
  \|u\|_{L^p L^q (\Rm \times \Rm^d)} \lesssim \|(u_0,u_1)\|_{\dot{H}^1 \times L^2}
 \]
 holds for all linear free wave $u = \mathbf{S}_L (u_0,u_1)$ with initial data $(u_0,u_1)$. Let us focus on the exterior region $\{(x,t): |x|>|t|\}$ and use the notation $\chi$ for the characteristic function of this region. We recall the fact that the inequality 
 \[
  \left(\sum_{j=-\infty}^\infty \|\chi_j u\|_{L^p L^q (\Rm \times \Rm^d)}^q\right)^{1/q} \leq \|\chi u\|_{L^p L^q (\Rm \times \Rm^d)} \leq \left(\sum_{j=-\infty}^\infty \|\chi_j u\|_{L^p L^q (\Rm \times \Rm^d)}^p\right)^{1/p}
 \]
 In the case $3 \leq d \leq 6$, we have $p \geq 2$. Therefore we always have 
 \[
   \|\chi u\|_{L^p L^q (\Rm \times \Rm^d)} \leq \left(\sum_{j=-\infty}^\infty \|\chi_j u\|_{L^p L^q (\Rm \times \Rm^d)}^2\right)^{1/2}.
 \]
 Thus our main theorem gives a stronger Strichartz estimates of radial free waves than the classic ones if we focus on the exterior region. If $d \geq 7$, the classic Strichartz estimate $\|u\|_{L^p L^q (\Rm \times \Rm^d)} \lesssim \|(u_0,u_1)\|_{\dot{H}^1 \times L^2}$ fails. In this case our main theorem gives a substitution of the classic ones. 
\end{remark}

\section{Strichartz estimates of radial free waves}

We first consider the case of free waves, i.e. we prove the main theorem when $F =0$.  Our main tool in this section is the radiation fields of free waves. 
\begin{theorem}[Radiation filed, see Duyckaerts et al. \cite{dkm3}] \label{radiation}
Assume that $d\geq 3$ and let $u$ be a solution to the free wave equation $\partial_t^2 u - \Delta u = 0$ with initial data $(u_0,u_1) \in \dot{H}^1 \times L^2(\Rm^d)$. Then 
\[
 \lim_{t\rightarrow \pm \infty} \int_{\Rm^d} \left(|\nabla u(x,t)|^2 - |u_r(x,t)|^2 + \frac{|u(x,t)|^2}{|x|^2}\right) dx = 0
\]
 and there exist two functions $G_\pm \in L^2(\Rm \times \mathbb{S}^{d-1})$ so that
\begin{align*}
 \lim_{t\rightarrow \pm\infty} \int_0^\infty \int_{\mathbb{S}^{d-1}} \left|r^{\frac{d-1}{2}} \partial_t u(r\theta, t) - G_\pm (r\mp t, \theta)\right|^2 d\theta dr &= 0;\\
 \lim_{t\rightarrow \pm\infty} \int_0^\infty \int_{\mathbb{S}^{d-1}} \left|r^{\frac{d-1}{2}} \partial_r u(r\theta, t) \pm G_\pm (r\mp t, \theta)\right|^2 d\theta dr & = 0.
\end{align*}
In addition, the maps $(u_0,u_1) \rightarrow \sqrt{2} G_\pm$ are a bijective isometries form $\dot{H}^1 \times L^2(\Rm^d)$ to $L^2 (\Rm \times \mathbb{S}^{d-1})$.
\end{theorem}
We usually call the functions $G_\pm$ radiation profiles of the corresponding free wave, or equivalently, of the initial data. It is clear that the free wave is radial if and only if $G_\pm$ is independent of the angle $\theta$.  We can give an explicit formula of free waves $u$ in term of their radiation profiles $G_-$ as below
 \[
  u(x,t) = \frac{1}{(2\pi)^\mu} \int_{\mathbb{S}^{d-1}} G_-^{(\mu-1)} (x \cdot \omega + t, \omega) d\omega,\qquad \mu = \frac{d-1}{2}.
 \]
The proof  can be found in \cite{shenradiation}. If the free wave is radial, or equivalently speaking, $G_-$ is independent of the angle, we may write $\omega=(\omega_1, \omega_2, \cdots, \omega_d) \in \mathbb{S}^{d-1}$ and obtain 
 \begin{align*}
 u(r e_1,t) =  (2\pi)^{-\mu} \int_{\mathbb{S}^{d-1}} G_-^{(\mu-1)} (r \omega_1+t) d\omega = \frac{\sigma_{d-2}}{(2\pi)^{\mu}}  \int_{-1}^{1} G_-^{(\mu-1)} (r \omega_1+t) (1-\omega_1^2)^{\mu-1} d\omega_1.
\end{align*}
Here $e_1=(1,0,\cdots,0)$ and $\sigma_{d-2}$ is the area of the sphere $\mathbb{S}^{d-2}$. We then integrate by parts
\begin{align*}
 u (r e_1,t)  = \frac{(-1)^{\mu-1} \sigma_{d-2}}{(2\pi)^{\mu} r^{\mu-1}} \int_{-1}^{1} G_-(r\omega_1+t) \left[\frac{d^{\mu-1}}{d\omega_1^{\mu-1}} (1-\omega_1^2)^{\mu-1}\right] d\omega_1.
\end{align*}
Letting $P_\frac{d-3}{2}$ be the Legendre polynomial of degree $\mu-1=(d-3)/2$ defined by
\[
 P_\frac{d-3}{2}(z) =  \frac{(-1)^{\mu-1} \sigma_{d-2}}{(2\pi)^{\mu}} \frac{d^{\mu-1}}{dz^{\mu-1}} (1-z^2)^{\mu-1} = \frac{1}{2^{\mu-1}(\mu-1)!} \frac{d^{\mu-1}}{dz^{\mu-1}} (z^2-1)^{\mu-1}
\]
and applying the change of variable $s=r\omega_1 + t$, we obtain
\begin{equation} \label{radial explicit formula}
 u(r,t) = r^{-(d-1)/2} \int_{-r+t}^{r+t} G_-(s) P_\frac{d-3}{2}\left(\frac{s-t}{r}\right) ds. 
\end{equation} 
We may also generalize the theory of radiation fields to free waves with initial data in the Sobolev spaces $\dot{H}^s \times \dot{H}^{s-1} (\Rm^d)$. 

\begin{proposition}\label{prop isometry}
 Given $\beta \in \Rm$, the maps from initial data to the radiation profiles $(u_0,u_1)\rightarrow \sqrt{2} G_\pm$ are isometries from the Hilbert space $\dot{H}^{\beta} (\Rm^d) \times \dot{H}^{\beta -1}(\Rm^d)$ to the space $L^2(\mathbb{S}^{d-1}, \dot{H}^{\beta-1}(\Rm))$. 
\end{proposition}
The proof of this proposition can be found in \cite{radiationsub}. In fact this is a direct consequence of the explicit formula for radiation fields in term of Fourier transform given in \cite{newradiation}.

\paragraph{Localized radiation profiles} We first assume that the radiation profile $G_-(s)$ is supported in $J_k = [-2^{k+1}, -2^k] \cup [2^k, 2^{k+1}]$ and consider the upper bound of $\|\chi_{j} u\|_{L^p L^q (\Rm^d \times \Rm)}$. We let $\chi_{k,r,t}$ be the characteristic function of the set $J_k \cap [-r+t,r+t]$ and rewrite the formula \eqref{radial explicit formula} as 
\[
 u(r,t) = r^{-(d-1)/2} \int_{\Rm} \chi_{k,r,t}(s) G_-(s) P_\frac{d-3}{2}\left(\frac{s-t}{r}\right) ds.
\]
This immediately gives 
\[
 |u(r,t)| \leq r^{-(d-1)/2} \|G_-\|_{\dot{H}^{\beta-1} (\Rm)} \left\|\chi_{k,r,t} (\cdot)P_\frac{d-3}{2}\left(\frac{\cdot-t}{r}\right) \right\|_{\dot{H}^{1-\beta}(\Rm)}.
\]
In order to give a reasonable $\dot{H}^{1-\beta}(\Rm)$ upper bound, we need to use the following lemma
\begin{lemma} \label{technical Hgamma}
 Let $g \in C^1(J)$ be a function defined in an interval $J = [a,b]$ and $\gamma \in (-1/2,1/2)$. We may also view $g$ as a function defined in $\Rm$ by defining $g(s) = 0$ for $s \notin J$. Then we have
 \begin{itemize}
  \item[(i)] If $\gamma \in (-1/2, 0]$, then $\|g\|_{\dot{H}^\gamma(\Rm)} \lesssim_\gamma |J|^{1/2-\gamma} \sup |g|$;
  \item[(ii)] If $\gamma \in (0,1/2)$, then $\|g\|_{\dot{H}^\gamma(\Rm)} \lesssim_\gamma |J|^{1/2-\gamma} (\sup |g| + |J|\sup |g'|)$.
 \end{itemize}
\end{lemma}
\begin{proof}
 We start by proving (i). Let $q\in (1,2]$ so that $1/q = 1/2 - \gamma$. A basic calculation shows that 
 \[
  \|g\|_{L^q (\Rm)} \lesssim_\gamma |J|^{1/q} \sup |g|.
 \]
 A Sobolev embedding $L^q(\Rm) \hookrightarrow \dot{H}^\gamma(\Rm)$ then finishes the proof of part (i). In order to verify (ii) we consider the Fourier transform
 \[
  \hat{g} (\xi) = c \int_a^b e^{is\xi} g(s) ds = \frac{ci e^{ia\xi}}{\xi} g(a) - \frac{ci e^{ib\xi}}{\xi} g(b) + \frac{ci}{\xi} \int_a^b e^{is \xi} g'(s) ds. 
 \] 
 Therefore we have 
 \[
  |\hat{g} (\xi)| \lesssim |\xi|^{-1} (\sup |g| + |J| \sup |g'|).
 \]
 This immediately gives that for any $R > 0$, we have 
 \begin{align*}
  \|g\|_{\dot{H}^\gamma}^2 & \lesssim R^{2\gamma} \int_{-R}^R |\hat{g}(\xi)|^2 d\xi + \int_{|\xi|>R} |\xi|^{2\gamma} |\hat{g}(\xi)|^2 d\xi\\
  & \lesssim R^{2\gamma} \|g\|_{L^2}^2 + (\sup |g| + |J| \sup |g'|)^2 \int_{|\xi|>R} |\xi|^{2\gamma-2} d\xi\\
  & \lesssim_\gamma R^{2\gamma} |J| (\sup |g|)^2 + R^{2\gamma-1} (\sup |g| + |J| \sup |g'|)^2.
 \end{align*}
 Finally we choose $R = |J|^{-1}$ and finishes the proof of (ii). 
\end{proof}
We next apply Lemma \ref{technical Hgamma} and obtain 
\[
  \left\|\chi_{k,r,t} (\cdot)P_\frac{d-3}{2}\left(\frac{\cdot-t}{r}\right) \right\|_{\dot{H}^{-\beta}(\Rm)} \lesssim_\beta (\min\{r, 2^k\})^{\beta-1/2}. 
\]
As a result we have the upper bound 
\[
 |u(r,t)| \lesssim_\beta r^{-(d-1)/2} (\min\{r, 2^k\})^{\beta-1/2} \|G_-\|_{\dot{H}^{\beta-1}}.
\]
We consider two cases: If $j\geq k-1$, then a straight-forward calculation gives (we recall the assumption $1/p + d/q = d/2 - \beta$)
\[
 \|\chi_j u\|_{L^p L^q (\Rm^d \times \Rm)} \lesssim (2^k/2^j)^{\beta-1/2} \|G_-\|_{\dot{H}^{\beta-1}}. 
\]
If $j \leq k-2$, then the explicit formula \eqref{radial explicit formula} implies that 
\[
 u(x,t) = 0, \qquad  \forall (x,t) \in \Omega_j \cap \{(x,t): |t| < 2^{k-2}\}. 
\]
We then conduct a straight-forward calculation and obtain 
\begin{align*}
 \|\chi_j u\|_{L^p L^q (\Rm^d \times \Rm)} & \lesssim (2^k)^{\beta-1/2} \|G_-\|_{\dot{H}^{\beta-1}} \left(\int_{|t|>2^{k-2}} \left(\int_{|t|+2^j}^{|t|+2^{j+1}} r^{-\frac{(d-1)q}{2}} r^{d-1} dr \right)^{p/q} dt\right)^{1/p} \\
 & \lesssim (2^j/2^k)^{1/q} \|G_-\|_{\dot{H}^{\beta-1}}.
\end{align*}
In summary we have that for any $j, k \in \mathbb{Z}$, if the radiation profile $G_-$ of a radial free wave $u$ is supported in $[-2^{k+1}, -2^k] \cup [2^k, 2^{k+1}]$, then 
\begin{equation} \label{single channel localized profile}
 \|\chi_j u\|_{L^p L^q (\Rm^d \times \Rm)} \lesssim_{p,q} c_{j-k} \|G_-\|_{\dot{H}^{\beta-1}}. 
\end{equation}
Here $\{c_n\}_{n\in \mathbb{Z}}$ is an $l^1$ sequence defined by 
\[
 c_n = \left\{\begin{array}{ll} 2^{(1/2-\beta)n}, & n \geq -1; \\ 2^{n/q}, & n \leq -2. \end{array}\right.
\]

\paragraph{Decomposition of profiles} Let $u$ be a radial free wave with initial data $(u_0,u_1)\in \dot{H}^\beta \times \dot{H}^{\beta-1}(\Rm^d)$. It immediately follows that its radiation profile $G_-(s)$ is a $\dot{H}^{\beta-1} (\Rm)$ function. We first introduce a technical lemma, whose proof is postponed to the Appendix. 

\begin{lemma} \label{decomposition of Hs}
 Let $\gamma \in (-1/2,1/2)$ and $G \in \dot{H}^\gamma (\Rm)$. If we define $G_k \in \dot{H}^\gamma (\Rm)$ by 
 \[
  G_k (s) = \left\{\begin{array}{ll} G(s), & \hbox{if} \; s \in [-2^{k+1}, -2^k] \cup [2^k, 2^{k+1}]; \\ 0, & \hbox{otherwise}; \end{array}\right.
 \]
 then we have the following series holds in the sense of $\dot{H}^\gamma (\Rm)$ strong convergence
 \[
  G = \sum_{k= - \infty}^\infty G_k(s). 
 \]
 In addition the following limit holds 
 \[
  \sum_{k=-\infty}^\infty \|G_k\|_{\dot{H}^\gamma (\Rm)}^2 \lesssim_\gamma \|G\|_{\dot{H}^\gamma (\Rm)}^2. 
 \]
\end{lemma}

\noindent We then decompose $G_-$ into the sum of $G_-^k$'s accordingly and let the corresponding radial free wave be $u_k$. We recall the upper bound given in \eqref{single channel localized profile} and obtain 
\[
 \|\chi_j u_k\|_{L^p L^q (\Rm^d \times \Rm)} \lesssim_{p,q} c_{j-k} \|G_-^k\|_{\dot{H}^{\beta-1}}.
\]
Therefore we may take a sum 
\[
 \|\chi_j u\|_{L^p L^q (\Rm^d \times \Rm)} \lesssim_{p,q} \sum_{k=-\infty}^\infty c_{j-k} \|G_-^k\|_{\dot{H}^{\beta-1}}.
\]
The right hand side can be viewed as a convolution of sequences $\{c_n\}_{n \in \mathbb{Z}}$ and $\{\|G_-^k\|_{\dot{H}^{\beta-1}}\}_{k\in \mathbb{Z}}$. Therefore we may apply Young's inequality and conclude 
\begin{align*}
 \left(\sum_{j=-\infty}^\infty \|\chi_j u\|_{L^p L^q (\Rm^d \times \Rm)}^2 \right)^{1/2} & \lesssim_{p,q} \|c_n\|_{l^1} \left(\sum_{k=-\infty}^\infty \|G_k\|_{\dot{H}^{\beta-1} (\Rm)}^2\right)^{1/2} \\
 & \lesssim_{p,q} \|G_-\|_{\dot{H}^{\beta-1} (\Rm)} \lesssim_{p,q} \|(u_0,u_1)\|_{\dot{H}^\beta \times \dot{H}^{\beta-1}(\Rm^d)}. 
\end{align*}

\section{Estimates of Inhomogeneous Term} 

Now we consider the radial solution $u$ to the wave equation $\partial_t^2 u - \Delta u = F$ with zero initial data. We first make the channel-localized decomposition 
\begin{align*}
  & F = \sum_{k = -\infty}^\infty F_k, \qquad |x|>|t|; & &F_k  = \chi_k F.&
\end{align*}
and define $u_k$ accordingly 
\[
 u_k = \int_{0}^t \mathbf{S}_L(t-\tau) (0, F_k(\cdot,\tau)) d\tau. 
\]
Then finite speed of propagation gives the decomposition 
\[
 \chi_j u = \sum_{k \geq j} \chi_j u_k. 
\]
We first find a upper bound of $\|\chi_j u_k\|_{L^p L^q}$ in term of $\|F_k\|_{L^{\tilde{p}} L^{\tilde{q}}}$. The idea is to find an explicit formula of $u_k$ in term of $F_k$. We start by considering the linear wave propagation operator $\mathbf{S}_L (t) (0,u_1)$. Here we assume $u_1$ is supported in the sphere shell $\{x: a<|x|<b\}$. We recall the explicit formula 
\[
 \mathbf{S}_L (t) (0,u_1) (x,t) = c_d \cdot  \left(\frac{1}{t}\frac{\partial}{\partial t}\right)^{(d-3)/2}\left(t^{d-2} \int_{\mathbb{S}^{d-1}} u_1(x + t\omega) d\omega\right).
\]
We let $r = |x+t\omega|$, use the radial assumption and apply a change of variable $\omega \rightarrow r$
\begin{align*}
 \int_{\mathbb{S}^{d-1}} u_1(x + t\omega) d\omega & = c_d \int_{|x|-t}^{|x|+t} u_1(r) \left[1 - \left(\frac{|x|^2 + t^2 - r^2}{2|x| t}\right)^2\right]^\frac{d-3}{2} \frac{rdr}{|x|t} \\
 & = c_d \frac{1}{|x|^{d-2} t^{d-2}} \int_{|x|-t}^{|x|+t} r u_1(r) \left[(|x|+t)^2 - r^2\right]^{\frac{d-3}{2}} \left[r^2 - (|x|-t)^2\right]^{\frac{d-3}{2}} dr.
\end{align*}
We plug this integral in the formula of $\mathbf{S}_L (t) (0,u_1)$ given above, and obtain
\begin{align*}
 \mathbf{S}_L (t) (0,u_1) (x,t) & = \frac{1}{|x|^{d-2}}\left(\frac{1}{t}\frac{\partial}{\partial t}\right)^{(d-3)/2} \left(\int_{|x|-t}^{|x|+t} r u_1(r) P_\frac{d-3}{2} (|x|,t,r) dr\right)\\
 & = \frac{1}{|x|^{d-2}} \sum_{m = 0}^{(d-3)/2} \frac{A_{d,m}}{t^{d-3-m}} \int_{|x|-t}^{|x|+t} r u_1(r) \left[\left(\frac{\partial}{\partial t}\right)^m P_\frac{d-3}{2} (|x|,t,r)\right]dr.
\end{align*}
Here $A_{d,m}$ are constants and $P_{\frac{d-3}{2}} (|x|,t,r)$ is a symmetric polynomial of three variables 
\[
 P_{\frac{d-3}{2}} (|x|,t,r) = \left[(|x|+t+r)(|x|+t-r)(r+|x|-t)(r+t-|x|)\right]^{\frac{d-3}{2}}. 
\]
We now give a point-wise estimate of $\mathbf{S}_L (t) (0,u_1)$. 

\begin{lemma} \label{pointwise L tilde q}
 Assume $u_1 \in L^{\tilde{q}}(\Rm^d)$ is a radial function supported in the sphere shell $\{x\in \Rm^d: a<|x|<b\}$. Here the positive constants $a<b$ satisfy $b/a \leq 2$. Then the free wave $u = \mathbf{S}_L (t) (0,u_1)$ satisfies 
\[
 |u(x,t)| \lesssim_d \frac{a^{(d-1)(1/2-1/\tilde{q})}}{|x|^{(d-1)/2}} (b-a)^{1-1/\tilde{q}} \|u_1\|_{L^{\tilde{q}}(\Rm^d)}, \qquad |x|>|t|>0. 
\]
\end{lemma}
\begin{proof}
The proof follows a straight-forward calculation. By symmetry we only need to deal with the positive time direction. We consider two cases, i.e. when time $t<a/4$ is small and when $t\geq a/4$ is large. Let us first consider the first case. By finite speed of propagation it suffice to consider the case $3a/4 \leq |x| \leq b+a/4$. We observe that if $r \in (|x|-t, |x|+t)$, then the inequalities $|t+r-|x||, ||x|+t-r|\leq 2t$ and $3a/2 \leq t+r+|x|, |x|+r-t \leq 5a$ hold. As a result we have 
 \[
  \left|\left(\frac{\partial}{\partial t}\right)^m P_\frac{d-3}{2} (|x|,t,r)\right| \lesssim_d t^{d-3-m} a^{d-3}.
 \]
 Plugging this upper bound in the explicit formula given above, we obtain 
 \begin{align*}
  |u(x,t)| & \lesssim_d \frac{1}{a} \int_{|x|-t}^{|x|+t} r |u_1(r)| dr \lesssim_d \frac{1}{a} \left(\int_{a}^{b} r^{\frac{\tilde{q} +1 -d}{\tilde{q}-1}} dr\right)^{1-1/{\tilde{q}}}\left(\int_{a}^{b} r^{d-1} |u_1(r)|^{\tilde{q}} dr\right)^{1/{\tilde{q}}}\\
  & \lesssim_d (b-a)^{1-1/\tilde{q}} a^{(1-d)/\tilde{q}} \|u_1\|_{L^{\tilde{q}}}\\
  & \lesssim_d \frac{a^{(d-1)(1/2-1/\tilde{q})}}{|x|^{(d-1)/2}} (b-a)^{1-1/\tilde{q}} \|u_1\|_{L^{\tilde{q}}(\Rm^d)}. 
 \end{align*}
 Next we consider the second case $t \geq a/4$. By finite speed of propagation, we always have $u(x,t) = 0$ if $|x|>b+t$. Therefore it suffices to consider the case $t<|x|<t+b$. In this case we have that if $r \in [|x|-t, |x|+t]\cap [a,b]$, then
 \begin{align*}
  &|x| \simeq_1 t \simeq_1 t+r+|x|;& &|x|+t -r \lesssim_1 |x|;& &t+r-|x|, |x|+r-t \lesssim_1 a.&
 \end{align*} 
 This implies that 
 \[
  \left|\left(\frac{\partial}{\partial t}\right)^m P_\frac{d-3}{2} (|x|,t,r)\right| \lesssim_d |x|^{d-3} a^{d-3-m}.
 \]
 Therefore we have 
  \begin{align*}
  |u(x,t)| & \lesssim_d  \frac{1}{|x|^{d-2}} \sum_{m = 0}^{(d-3)/2} \frac{|A_{d,m}|}{|x|^{d-3-m}} \int_{|x|-t}^{|x|+t} r |u_1(r)| |x|^{d-3} a^{d-3-m}dr\\
 &\lesssim_d \frac{a^{(d-3)/2}}{|x|^{(d-1)/2}} \int_{a}^{b} r |u_1(r)| dr\\
 & \lesssim_d \frac{a^{(d-3)/2}}{|x|^{(d-1)/2}} \left(\int_{a}^{b} r^{\frac{\tilde{q} +1 -d}{\tilde{q}-1}} dr\right)^{1-1/{\tilde{q}}}\left(\int_{a}^{b} r^{d-1} |u_1(r)|^{\tilde{q}} dr\right)^{1/{\tilde{q}}}\\
  & \lesssim_d \frac{a^{(d-1)(1/2-1/\tilde{q})}}{|x|^{(d-1)/2}} (b-a)^{1-1/\tilde{q}} \|u_1\|_{L^{\tilde{q}}}.
  \end{align*}
\end{proof}

\paragraph{$L^{\tilde{q}}-L^q$ estimate} Now we assume $F_k \in L^{\tilde{p}} L^{\tilde{q}} (\Rm \times \Rm^d)$ is supported in $\Omega_k$. If $j \leq k$, then we have 
\[
 \chi_j u_k = \int_{0}^t \chi_j \mathbf{S}_L(t-\tau) (0, F_k(\cdot,\tau)) d\tau. 
\]
By Lemma \ref{pointwise L tilde q} we have 
\[
 \left\|\chi_j \mathbf{S}_L(t-\tau) (0, F_k(\cdot,\tau))\right\|_{L^q (\Rm^d)} \lesssim_d \frac{(2^k + \tau)^{(d-1)(1/2-1/\tilde{q})}}{(2^j +t)^{(d-1)(1/2-1/q)}} (2^k)^{1-1/\tilde{q}} (2^j)^{1/q}\|F_k(\cdot,\tau)\|_{L^{\tilde{q}}(\Rm^d)}.
\]
We may combine this estimate with the finite speed of propagation and obtain
\begin{align*}
 \|\chi_j u_k(\cdot,t)\|_{L^q(\Rm^d)} &\lesssim_d \int_0^{\tau_+(j,k,t)} \frac{(2^k + \tau)^{(d-1)(1/2-1/\tilde{q})}}{(2^j +t)^{(d-1)(1/2-1/q)}} (2^k)^{1-1/\tilde{q}} (2^j)^{1/q}\|F_k(\cdot,\tau)\|_{L^{\tilde{q}}(\Rm^d)} d\tau\\
 & \lesssim_d \frac{(2^k)^{1-1/\tilde{q}} (2^j)^{1/q}}{(2^j +t)^{(d-1)(1/2-1/q)}} \left\|(2^k + \tau)^{(d-1)(1/2-1/\tilde{q})}\right\|_{L_\tau^{\tilde{p}'}([0,\tau_+(j,k,t)])} \|F_k\|_{L^{\tilde{p}} L^{\tilde{q}}}.
\end{align*}
Here $\tilde{p}'$ is the conjugate index of $\tilde{p}$ and $\tau_+(j,k,t)$ is defined by 
\[
 \tau_+(j,k,t) = \left\{\begin{array}{ll} t, & \hbox{if}\; j \in \{k-1,k\}; \\ t+2^{j} - 2^{k-1}, & \hbox{if}\; j \leq k-2, \; t > 2^{k-1} - 2^{j}; \\ 0, & \hbox{if}\; j \leq k-2, \; t \leq 2^{k-1} - 2^{j}. \end{array} \right.
\]
Except for the trivial case $\tau_+ = 0$ we always have 
\[
 \left\|(2^k + \tau)^{(d-1)(1/2-1/\tilde{q})}\right\|_{L_\tau^{\tilde{p}'}([0,\tau_+(j,k,t))} \leq \left\|y^{(d-1)(1/2-1/\tilde{q})}\right\|_{L_y^{\tilde{p}'}([2^k,t+2^{j}+2^{k-1}])}.
\]
The upper bound of the right hand side depends on the value of $(d-1)(1/2-1/\tilde{q})$ and $\tilde{p}'$. There are three cases. 

\paragraph{Case 1} If $(d-1)(1/2-1/\tilde{q}) + 1/\tilde{p}' < 0$, then we have  
\[
 \left\|y^{(d-1)(1/2-1/\tilde{q})}\right\|_{L_y^{\tilde{p}'}([2^k,t+2^{j}+2^{k-1}])} \lesssim (2^k)^{(d-1)(1/2-1/\tilde{q}) + 1/\tilde{p}'}. 
\]
Thus we have 
\[
  \|\chi_j u_k(\cdot,t)\|_{L^q(\Rm^d)}  \lesssim \frac{(2^k)^{\beta-1/2} (2^j)^{1/q}}{(2^j +t)^{(d-1)(1/2-1/q)}}  \|F_k\|_{L^{\tilde{p}} L^{\tilde{q}}}.
\]
Here we utilize the assumption $1/\tilde{p} + d/\tilde{q} = d/2 + 2 -\beta$. A direct calculation (as well as finite speed of propagation) shows
\begin{align*}
 \|\chi_j u_k\|_{L^p L^q(\Rm^+ \times \Rm^d)} & \lesssim (2^k)^{\beta-1/2} (2^j)^{1/q} \left\|(2^j +t)^{-(d-1)(1/2-1/q)}\right\|_{L^p ([t_-(j,k), +\infty))} \|F_k\|_{L^{\tilde{p}} L^{\tilde{q}}} \\
 & \lesssim (2^k)^{\beta-1/2} (2^j)^{1/q} (2^k)^{-(d-1)(1/2-1/q) + 1/p} \|F_k\|_{L^{\tilde{p}} L^{\tilde{q}}} \\
 & \lesssim 2^{(j-k)/q} \|F_k\|_{L^{\tilde{p}} L^{\tilde{q}}}. 
\end{align*}
Here the lower limit 
\[
 t_-(j,k) = \left\{\begin{array}{ll} 0, & j \in \{k,k-1\}; \\ 2^{k-1} - 2^{j}, & j \leq k-2; \end{array}\right.
\]
guarantees that $2^j + t_-(j,k) \simeq 2^k$ and our assumptions $1/p + d/q = d/2 -\beta$, $p,q \in [1,+\infty)$ implies 
\[
 -(d-1)(1/2-1/q) + 1/p = 1/2 - \beta - 1/q < 0. 
\]

\paragraph{Case 2} If $(d-1)(1/2-1/\tilde{q}) + 1/\tilde{p}' > 0$, then we have (please note that unless $\tau_+(j,k,t) = 0$ we always have $t+2^j > 2^{k-1}$)
\[
 \left\|y^{(d-1)(1/2-1/\tilde{q})}\right\|_{L_y^{\tilde{p}'}([2^k,t+2^{j}+2^{k-1}])} \lesssim (t+2^j)^{(d-1)(1/2-1/\tilde{q}) + 1/\tilde{p}'} = (t+2^j)^{1/\tilde{q} +\beta -3/2}. 
\]
Thus we have
\[
  \|\chi_j u_k(\cdot,t)\|_{L^q(\Rm^d)}  \lesssim  (2^k)^{1-1/\tilde{q}} (2^j)^{1/q}(2^j +t)^{-(d-1)(1/2-1/q)+1/\tilde{q} +\beta -3/2}  \|F_k\|_{L^{\tilde{p}} L^{\tilde{q}}}
\]
A similar argument to case 1 gives 
\begin{align*}
 \|\chi_j u_k\|_{L^p L^q} & \lesssim (2^k)^{1-1/\tilde{q}} (2^j)^{1/q} \left\|(2^j +t)^{-(d-1)(1/2-1/q)+1/\tilde{q} +\beta -3/2}\right\|_{L^p ([t_-(j,k), +\infty))} \|F_k\|_{L^{\tilde{p}} L^{\tilde{q}}} \\
 & \lesssim (2^k)^{1-1/\tilde{q}} (2^j)^{1/q} (2^k)^{-(d-1)(1/2-1/q) + 1/p +1/\tilde{q} +\beta -3/2} \|F_k\|_{L^{\tilde{p}} L^{\tilde{q}}} \\
 & \lesssim 2^{(j-k)/q} \|F_k\|_{L^{\tilde{p}} L^{\tilde{q}}}. 
\end{align*}
Here we use the fact 
\[
 -(d-1)(1/2-1/q) + 1/p +1/\tilde{q} +\beta -3/2 = 1/\tilde{q} -1/q - 1 < 0.  
\]

\paragraph{Case 3} Finally if $(d-1)(1/2-1/\tilde{q}) + 1/\tilde{p}' = 0$, then we have 
\[
 \left\|y^{(d-1)(1/2-1/\tilde{q})}\right\|_{L_y^{\tilde{p}'}([2^k,t+2^{j}+2^{k-1}])} \lesssim \left(\ln \frac{t+2^{j}+2^{k-1}}{2^k}\right)^{1-1/\tilde{p}}. 
\]
Thus (again we recall that $t+2^j > 2^{k-1}$ unless $\tau_+(j,k,t) = 0$)
\[
 \|\chi_j u_k(\cdot,t)\|_{L^q(\Rm^d)} \lesssim \frac{(2^k)^{1-1/\tilde{q}} (2^j)^{1/q}}{(t+2^j +2^{k-1})^{(d-1)(1/2-1/q)}} \left|\ln \frac{t+2^{j}+2^{k-1}}{2^k}\right|^{1-1/\tilde{p}} \|F_k\|_{L^{\tilde{p}} L^{\tilde{q}}}
\]
Again we have 
\begin{align*}
 \|\chi_j u_k\|_{L^p L^q} & \lesssim (2^k)^{1-1/\tilde{q}} (2^j)^{1/q} \left\|\tilde{t}^{-(d-1)(\frac{1}{2}-\frac{1}{q})} \left(\ln \frac{\tilde{t}}{2^k}\right)^{1-\frac{1}{\tilde{p}}} \right\|_{L^p ([t_-(j,k)+2^j+2^{k-1}, +\infty))} \|F_k\|_{L^{\tilde{p}} L^{\tilde{q}}} \\
 & \lesssim (2^k)^{1-1/\tilde{q}} (2^j)^{1/q} \left\|\tilde{t}^{-(d-1)(\frac{1}{2}-\frac{1}{q})} \left(\ln \frac{\tilde{t}}{2^k}\right)^{1-\frac{1}{\tilde{p}}} \right\|_{L^p ([2^k, +\infty))} \|F_k\|_{L^{\tilde{p}} L^{\tilde{q}}} \\
 & \lesssim (2^k)^{1-1/\tilde{q}} (2^j)^{1/q} (2^k)^{-(d-1)(1/2-1/q) + 1/p} \|F_k\|_{L^{\tilde{p}} L^{\tilde{q}}} \\
 & \lesssim 2^{(j-k)/q} \|F_k\|_{L^{\tilde{p}} L^{\tilde{q}}}. 
\end{align*}
Here we recall $-(d-1)(1/2-1/q) + 1/p = 1/2 -\beta - 1/q < 0$ and observe that our assumptions $(d-1)(1/2-1/\tilde{q}) + 1/\tilde{p}' = 0$ and $1/\tilde{p} +d/\tilde{q} = d/2 -\beta +2$ imply $1-1/\tilde{q} = \beta - 1/2$. 

\paragraph{Summary} Combining three cases, we have for any $j \leq k$,
\[
 \|\chi_j u_k\|_{L^p L^q(\Rm^+ \times \Rm^d)} \lesssim 2^{(j-k)/q} \|F_k\|_{L^{\tilde{p}} L^{\tilde{q}}(\Rm\times \Rm^d)}
\]
Here the explicit constant does not depend on $j,k$ or $F$. We take a sum and obtain 
\[
  \|\chi_j u_k\|_{L^p L^q(\Rm^+ \times \Rm^d)} \lesssim \sum_{k \geq j} 2^{(j-k)/q} \|F_k\|_{L^{\tilde{p}} L^{\tilde{q}}(\Rm\times \Rm^d)}. 
\]
The right hand side is actually the convolution of an $l^2$ sequence $\{\|F_k\|_{L^{\tilde{p}} L^{\tilde{q}}(\Rm\times \Rm^d)}\}_k$ and an $l^1$ sequence $\{c_n\}_{n \in \mathbb{Z}}$ defined by
\[
 c_n = \left\{\begin{array}{ll} 2^{n/q}, & n\leq 0; \\ 0, & n > 0. \end{array}\right.
\]
Finally we apply Young's inequality to conclude that 
\[
 \left(\sum_{j = -\infty}^\infty \|\chi_j u_k\|_{L^p L^q(\Rm^+ \times \Rm^d)}^2 \right)^{1/2} \lesssim \left(\sum_{j = -\infty}^\infty \|F_k\|_{L^{\tilde{p}} L^{\tilde{q}}(\Rm\times \Rm^d)}^2 \right)^{1/2}. 
\]
The negative time direction is similar. 

\section{Application}

Now we consider the exterior solutions to the wave equation 
\begin{equation} \label{cp1}
 \left\{\begin{array}{ll} \partial_t^2 u - \Delta u = F(t,x,u), & (x,t) \in \Rm^d \times \Rm; \\ (u,u_t)|_{t=0} = (u_0,u_1) \in \dot{H}^1 \times L^2(\Rm^d); & \end{array}\right. 
\end{equation} 
with radial data. Here the nonlinear term $F$ satisfies 
\begin{itemize}
 \item $F$ is a radial function of $x$; 
 \item The following inequalities hold
 \begin{align*}
  &|F(t,x,u)| \leq C |u|^\frac{d+2}{d-2};& &|F(t,x,u)-F(t,x,v)| \leq C(|u|^\frac{4}{d-2} + |v|^\frac{4}{d-2}) |u-v|.&
 \end{align*}
\end{itemize}
Let $I$ be a time interval containing zero and $R\geq 0$ be a constant. We say that a function $u$ defined in the exterior region $\Omega_R = \{(x,t) \in \Rm^d \times I: |x|>R+|t|\}$ is an exterior solution to the equation given above if and only if 
\[
 u(x,t) = \mathbf{S}_L(u_0,u_1) + \int_0^t \frac{\sin (t-\tau)\sqrt{-\Delta}}{\sqrt{-\Delta}} [\chi_R F(\tau,x,u)] d\tau, \qquad |x|>R+|t|.
\]
Here we add the characteristic function $\chi_R$ of the exterior region $\Omega_R$ to emphasize that $F(\tau,x,u)$ is only defined in $\Omega_R$. For convenience we define the following norms for functions $u, F$ defined in $\Omega_R$
\begin{align*}
 \|u\|_{Y(I)} & = \left(\sum_{j=-\infty}^\infty \|\chi_j u\|_{L^\frac{d+2}{d-2} L^\frac{2(d+2)}{d-2} (I \times \Rm^d)}^2 \right)^{1/2}; \\
 \|F\|_{Z(I)} & = \left(\sum_{j=-\infty}^\infty \|\chi_j F\|_{L^1 L^2 (I \times \Rm^d)}^2 \right)^{1/2}. 
\end{align*}
Here the values of $u, F$ at a point outside $\Omega_R$ are viewed as zero. Our main theorem guarantees that the function $u$ defined by
\[
 u(x,t) = \mathbf{S}_L(u_0,u_1) + \int_0^t \frac{\sin (t-\tau)\sqrt{-\Delta}}{\sqrt{-\Delta}} (\chi_R F) d\tau, \qquad |x|>R+|t|;
\]
satisfies the Strichartz estimate
\begin{equation}
 \|u\|_{Y(I)} \leq \|\chi_R \mathbf{S}_L(u_0,u_1)\|_{Y(I)} + C \|F\|_{Z(I)}^{\frac{d+2}{d-2}}. 
\end{equation} 
Here we have 
\[
 \|\chi_R \mathbf{S}_L(u_0,u_1)\|_{Y(I)} \leq \left(\sum_{j > \log_2 R - 1} \|\chi_j \mathbf{S}_L(u_0,u_1)\|\right)^{1/2} \lesssim \|(u_0,u_1)\|_{\dot{H}^1 \times L^2}. 
\]
In addition, a straight-forward calculation shows that ($p = \frac{d+2}{d-2}$, $q = \frac{2(d+2)}{d-2}$)
\begin{align*}
 \|\chi_j F(t,x,u)- \chi_j F(t,x,v)\|_{L^1 L^2} & \leq C \left\|\chi_j (|u|^{4/(d-2)}+|v|^{4/(d-2)})|u-v|\right\|_{L^1 L^2}\\
 & \leq C \left(\|\chi_j u\|_{L^p L^q}^{\frac{4}{d-2}}+ \|\chi_j v\|_{L^p L^q}^{\frac{4}{d-2}}\right) \|\chi_j (u-v)\|_{L^p L^q}
\end{align*}
Thus we have 
\begin{align*}
 \|F(t,x,u)-F(t,x,v)\|_{Z(I)} & \leq C \left(\sum_{j=-\infty}^\infty \|\chi_j u\|_{L^p L^q}^{2p} \right)^{\frac{p-1}{2p}} \left(\sum_{j=-\infty}^\infty \|\chi_j (u-v)\|_{L^p L^q}^{2p} \right)^{\frac{1}{2p}}\\
 & \qquad + C\left(\sum_{j=-\infty}^\infty \|\chi_j v\|_{L^p L^q}^{2p} \right)^{\frac{p-1}{2p}} \left(\sum_{j=-\infty}^\infty \|\chi_j (u-v)\|_{L^p L^q}^{2p} \right)^{\frac{1}{2p}}\\
 & \leq C\left(\|u\|_{Y(I)}^{p-1} + \|v\|_{Y(I)}^{p-1}\right) \|u-v\|_{Y(I)}.
\end{align*}
Letting $v=0$, we also have 
\[
 \|F(t,x,u)\|_{Z(I)} \leq C \|u\|_{Y(I)}^p. 
\]
Combining these inequalities with a standard fixed-point argument, we may obtain the following 
\begin{proposition}
 Assume that $d\geq 3$ is an odd integer. Let $(u_0,u_1)\in \dot{H}^1 \times L^2(\Rm^d)$ be radial initial data. 
 \begin{itemize} 
  \item There exists a constant $\delta>0$, so that if $\|\chi_R \mathbf{S}_L(u_0,u_1)\|_{Y(I)} < \delta$, then there exists a unique exterior solution $u$ defined in $\Omega_R \cap (\Rm^d \times I)$ satisfying $\|u\|_{Y(I)} < +\infty$. 
  \item In addition, if $v$ is another solution to (CP1) with initial data $(v_0,v_1)\in \dot{H}^1 \times L^2(\Rm^d)$ so that $\|\chi_R \mathbf{S}_L(u_0,u_1)\|_{Y(I)} < \delta$, then we have 
  \[
   \|u-v\|_{Y(I)} \leq 2 \|\chi_R \mathbf{S}_L (u_0-v_0, u_1-v_1)\|_{Y(I)}. 
  \]
  \item Given any radial initial data $(u_0,u_0) \in \dot{H}^1 \times L^2(\Rm^d)$ and $R\geq 0$, there exists a small time $T>0$, so that there exists a unique exterior solution defined in $\Omega_R \cap (\Rm^d \times [-T,T])$ satisfying $\|u\|_{Y([-T,T])} < +\infty$.
  \item Given any radial initial data $(u_0,u_0) \in \dot{H}^1 \times L^2(\Rm^d)$, there exists a radius $R \geq 0$, so that there exists a unique exterior solution defined in $\Omega_R$ satisfying $\|u\|_{Y(\Rm)} < +\infty$.
 \end{itemize}
\end{proposition}

\section*{Appendix}
In this appendix we give a proof of Lemma \ref{decomposition of Hs}. First of all, we recall that the cut-off operator $f \rightarrow \chi_{[0,+\infty)} f$ is a bounded operator from $\dot{H}^\beta (\Rm)$ to itself. Please see \cite{cutoff}, for instance. By applying the reflection and translation we obtain that 
\begin{equation} \label{cutoff operator}
 \|\chi_J f\|_{\dot{H}^\beta(\Rm)} \lesssim \|f\|_{\dot{H}^\beta(\Rm)}. 
\end{equation} 
Here $J$ is an arbitrary interval (finite or infintie) and the implicit constant does not depend on the choice of interval $J$. We first verify that 
\begin{equation} \label{decomposition in H beta}
  G = \sum_{k= - \infty}^{\infty} G_k(s). 
\end{equation}
By the universal boundedness of the cut-off operators given above, it suffices to consider a dense subset of the space $\dot{H}^\beta(\Rm)$. It is clear that the decomposition above holds for $G \in C_0^\infty(\Rm)$. This verifies \eqref{decomposition in H beta}. The majority of this section is devoted to the proof the inequality 
\[
 \sum_{k=-\infty}^{\infty} \|G_k\|_{\dot{H}^\beta(\Rm)}^2 \lesssim \|G\|_{\dot{H}^\beta(\Rm)}^2. 
\]
By symmetry it suffices to show that if $\beta \in (-1/2,1/2)$, then
\[
 \sum_{k = -\infty}^\infty \|\chi_{[2^k, 2^{k+1}]} f\|_{\dot{H}^\beta (\Rm)}^2 \lesssim \|f\|_{\dot{H}^\beta (\Rm)}^2. 
\]
Here we use that notation $\chi_{A}$ for the characteristic function of a set $A \subset \Rm$. By \eqref{cutoff operator} it suffice to show that 
\[
 \sum_{k = -\infty}^\infty \|\phi_k f\|_{\dot{H}^\beta (\Rm)}^2 \lesssim \|f\|_{\dot{H}^\beta (\Rm)}^2. 
\]
Here the functions are defined by 
\[
 \phi_k(x) = \left\{\begin{array}{ll} 1, & x\in [2^k, 2^{k+1}]; \\ 2^{1-k} x -1, & x \in [2^{k-1},2^k]; \\
 3-2^{-k} x, & x\in [2^{k+1}, 3\cdot 2^k]; \\ 0, & \hbox{otherwise}. \end{array}\right.
\]
Turning to the Fourier side (let $\hat{g} = f$), we only need to show that (please note that this clearly holds if $\beta = 0$.)
\begin{equation} \label{Fside}
 \sum_{k=-\infty}^\infty \int_{\Rm} |(\varphi_k \ast g)(\xi)|^2 |\xi|^{2\beta} d\xi \lesssim \int_\Rm |g(\xi)|^2 |\xi|^{2\beta} d\xi. 
\end{equation}
Here $g \in L^2(\Rm; |\xi|^{2\beta} d\xi)$ and $\varphi_k$ is the inverse Fourier transform of $\phi_k$. We may find two upper bounds of $|\varphi_k|$. On one hand, we have 
\[
 |\varphi_k(\xi)| \lesssim \left|\int_{2^{k-1}}^{3\cdot 2^k} \phi_k (x) e^{-ix\xi} dx \right| \lesssim 2^k.
\]
On the other hand, we may integrate by parts and obtain 
\begin{align*}
 |\varphi_k(\xi)| \lesssim \left|\frac{1}{\xi} \int_\Rm \phi'_k(x) e^{-ix\xi} dx \right| \lesssim 2^{-k} |\xi|^{-2}. 
 \end{align*}
 In summary, we always have 
 \[
   |\varphi_k(\xi)| \lesssim \min\{2^k, 2^{-k} |\xi|^{-2}\}. 
 \]
 We decompose the real line $\Rm = \cup I_m$ with $I_m = [-2^{m+1},-2^m] \cup [2^m, 2^{m+1}]$; and decompose $g$ accordingly 
 \begin{align*}
  &g = \sum_{n=-\infty}^\infty g_n;& &g_n(\xi)= \left\{\begin{array}{ll} g(\xi), & \xi \in I_n; \\ 0, & \xi\notin I_n. \end{array}\right.&
 \end{align*}
 We next write 
 \begin{align*}
  \sum_{k=-\infty}^\infty \int_{\Rm} |(\varphi_k \ast g)(\xi)|^2 |\xi|^{2\beta} d\xi & = \sum_{m = -\infty}^\infty \sum_{k=-\infty}^\infty \int_{I_m} \left|\left(\varphi_k \ast \sum_{n} g_n\right)(\xi)\right|^2 |\xi|^{2\beta} d\xi \lesssim J_1 + J_2. 
 \end{align*}
 Here $J_1$ and $J_2$ are defined by 
 \begin{align*}
  J_1 & = \sum_{m = -\infty}^\infty \sum_{k=-\infty}^\infty \int_{I_m} \left|\left(\varphi_k \ast \sum_{|n-m|\leq 1} g_n\right)(\xi)\right|^2 |\xi|^{2\beta} d\xi;\\
  J_2 & = \sum_{m = -\infty}^\infty \sum_{k=-\infty}^\infty \int_{I_m} \left|\left(\varphi_k \ast \sum_{|n-m|\geq 2} g_n\right)(\xi)\right|^2 |\xi|^{2\beta} d\xi.
 \end{align*}
 We first find the upper bound of $J_1$:
 \begin{align*}
  J_1 & \lesssim \sum_{\tau = -1,0,1} \left[\sum_{m = -\infty}^\infty \sum_{k=-\infty}^\infty \int_{I_m} 2^{2m\beta} \left|\left(\varphi_k \ast g_{m+\tau}\right)(\xi)\right|^2  d\xi\right] \\
  & \lesssim \sum_{\tau = -1,0,1} \left[\sum_{m = -\infty}^\infty \sum_{k=-\infty}^\infty \int_{\Rm} \left|\left(\varphi_k \ast 2^{m\beta} g_{m+\tau}\right)(\xi)\right|^2  d\xi\right] \\
  & \lesssim \sum_{\tau = -1,0,1}  \left[\sum_{m = -\infty}^\infty \int_{\Rm} \left|2^{m\beta} g_{m+\tau}(\xi)\right|^2 d\xi \right] \\
  & \lesssim \int_{\Rm} |g(\xi)|^2 |\xi|^{2\beta} d\xi. 
 \end{align*}
 Here we use the fact that \eqref{Fside} holds for $\beta=0$. Next we deal with $J_2$. We first observe that if $|m-n|\geq 2$, then the distance of $I_m$ and $I_n$ is about $2^{\max\{m,n\}}$. Thus if $\xi \in I_m$, then 
 \begin{align*}
  \left|(\varphi_k \ast g_n)(\xi)\right| & = \left|\int_{I_n} \varphi_k(\xi - \xi') g_n (\xi') d\xi' \right| \lesssim \int_{I_n} \min\{2^k, 2^{-k} |\xi - \xi'|^{-2}\} |g_n(\xi')| d\xi'\\
  & \lesssim \min\{2^k, 2^{-k-2m}, 2^{-k-2n}\} 2^{n/2} \|g_n\|_{L^2}\\
  & \lesssim 2^{-k+\frac{3}{2} \min\{k,-m\} - \Phi(\min\{k,-m\}+n)} \|g_n\|_{L^2}
 \end{align*}
 Here the function $\Phi$ is defined by 
 \[
  \Phi(y) = \left\{\begin{array}{ll} (-1/2) y, & y < 0; \\ (3/2) y, & y \geq 0. \end{array} \right.
 \]
 Therefore we have 
 \[
   \left|\left(\varphi_k \ast \sum_{|n-m|\geq 2} g_n\right)(\xi)\right| \lesssim \sum_{|n-m|\geq 2} 2^{-k+\frac{3}{2} \min\{k,-m\} - \Phi(\min\{k,-m\}+n)} \|g_n\|_{L^2}.
 \]
 We fix a constant $\gamma = 1^- > 2|\beta|$ and obtain 
 \[
 \left|\left(\varphi_k \ast \sum_{|n-m|\geq 2} g_n\right)(\xi)\right|^2 \lesssim \sum_{|n-m|\geq 2} 2^{-2k+3 \min\{k,-m\} - 2\gamma \Phi(\min\{k,-m\}+n)} \|g_n\|_{L^2}^2.
 \]
 We observe the fact that if we fix $m,n$, then the coefficient $2^{-2k+3 \min\{k,-m\} - 2\gamma \Phi(\min\{k,-m\}+n)}$ in the sum above takes its maximum value at $k = \min\{-m,-n\}$ and decays exponentially when $k$ moves away from $\min\{-m,-n\}$. Therefore we have 
 \[
  \sum_{k = -\infty}^\infty \left|\left(\varphi_k \ast \sum_{|n-m|\geq 2} g_n\right)(\xi)\right|^2 \lesssim \sum_{|n-m|\geq 2} 2^{(1+\gamma)\min\{-m,-n\} + \gamma n}\|g_n\|_{L^2}^2
 \]
 We then integrate and obtain 
 \[
  \sum_{k=-\infty}^\infty \int_{I_m} \left|\left(\varphi_k \ast \sum_{|n-m|\geq 2} g_n\right)(\xi)\right|^2 |\xi|^{2\beta} d\xi \lesssim \sum_{n=-\infty}^\infty 2^{(1+\gamma)\min\{-m,-n\} + \gamma n + m + 2\beta m}\|g_n\|_{L^2}^2
 \]
 Again we observe that if we fix $n$, then the coefficient $2^{(1+\gamma)\min\{-m,-n\} + \gamma n + m + 2\beta m}$ takes its maximum at $m=n$ and decays exponentially as $m$ moves away from $n$. Therefore we may take a sum in $m$ and obtain 
\[
 J_2 \lesssim \sum_{n=-\infty}^\infty 2^{2\beta n} \|g_n\|_{L^2}^2 \lesssim \int_\Rm |g(\xi)|^2 |\xi|^{2\beta} d\xi. 
\]
Finally we combine $J_1$ with $J_2$ to finish the proof.

\section*{Acknowledgement}
The second author is financially supported by National Natural Science Foundation of China Project 12071339.

\end{document}